\title{Compositions of random transpositions}
\author{Oded Schramm}
\documentclass[12pt]{article}
\usepackage{amsmath}
\usepackage{amsthm}
\usepackage{amsfonts}
\usepackage{graphicx}

\input labelfig.tex
\input epsf.tex

\IfFileExists{hyperref.sty}{\usepackage{hyperref}}{}

\newif\ifdraft
\drafttrue
\draftfalse
\def\note#1/{\ifdraft {\bf [#1]}\fi}
\numberwithin{equation}{section}
\numberwithin{figure}{section}

\newtheorem{theorem}{Theorem}
\numberwithin{theorem}{section}
\newtheorem{corollary}[theorem]{Corollary}
\newtheorem{lemma}[theorem]{Lemma}

\theoremstyle{definition}
\theoremstyle{definition}

\def\E{\mathbf E}
\def\P{\mathbf P}
\def\eref#1{(\ref{#1})}
\def\QED{\qed\medskip}

\newcommand{\R}{\mathbb{R}}

\newcommand{\Z}{\mathbb{Z}}
\newcommand{\N}{\mathbb{N}}
\newcommand{\Np}{{\mathbb{N}_+}}

\def\ceil#1{\lceil{#1}\rceil}
\def\bceil#1{\bigl\lceil{#1}\bigr\rceil}
\def\floor#1{\lfloor{#1}\rfloor}

\def\Ito/{It\^o}
\def \eps {\epsilon}
\def\md{\mid}
\def\Bb#1#2{{\def\md{\bigm| }#1\bigl[#2\bigr]}}
\def\BB#1#2{{\def\md{\Bigm| }#1\Bigl[#2\Bigr]}}
\def\Bs#1#2{{\def\md{\mid}#1[#2]}}
\def\Pb{\Bb\P}
\def\Eb{\Bb\E}
\def\PB{\BB\P}
\def\EB{\BB\E}
\def\Ps{\Bs\P}

\def\ev#1{{\mathcal{#1}}}
\def \proof {{ \medbreak \noindent {\bf Proof.} }}
\def\proofof#1{{ \medbreak \noindent {\bf Proof of #1.} }}

\def\TT{{q}}
\def\SS{{\mathfrak X}}
\def\cc{{\gamma}}

\def\Watterson{MR58:20594a}
\def\Dudley{MR91g:60001}
\def\DMZZ{math.PR/0305313}
\def\DiaconisShahshahani{MR82h:60024}
\def\BerestyckiDurrett{math.PR/0403259}

\def\Angel{AngelCyclicTime}
\def\Toth{MR94h:82062}
\def\Holst{HolstPD}
\def\SpencerLectures{MR95c:05113}
\def\AlonSpencer{MR2003f:60003}
\def\JLR{MR2001k:05180}
\def\riffle{MR1325269}

\ifx\hyperlink\undefined
\def\fff#1{&{{\pageref{#1}}}\cr}
\def\hfff#1{\label{#1}}
\else
\def\fff#1{&{\hyperlink{#1}{\pageref*{#1}}}\cr}
\def\hfff#1{\label{#1}\hypertarget{#1}}
\fi

\begin{document}
\maketitle

\begin{abstract}
Let $Y=(y_1,y_2,\dots)$, $y_1\ge y_2\ge\cdots$, be the list of sizes of the cycles
in the composition of $c\,n$ transpositions on the
set $\{1,2,\dots,n\}$.
We prove that if $c>1/2$ is constant and $n\to\infty$,  the distribution of $f(c)\,Y/n$
converges to $PD(1)$, the Poisson-Dirichlet distribution with paramenter $1$,
where the function $f$ is known explicitly.
A new proof is presented of the theorem by
Diaconis, Mayer-Wolf, Zeitouni and Zerner stating that the $PD(1)$ measure
is the unique invariant measure for the uniform coagulation-fragmentation process.
\end{abstract}

\section{Introduction}

Consider the composition \hfff{d.pi}{$\pi_t=T_t\circ T_{t-1}\circ\cdots T_2\circ T_1$}
of random, uniform, independent traspositions \hfff{d.T}{$T_j$} of 
\hfff{d.V}{$V:=\{1,2,\dots,n\}$}.
How large must $t$ be in order for $\pi_t$ to \lq\lq look like\rq\rq\ a random-uniform
permutation $\pi$ of $V$?
As we will see, the answer depends on the precise meaning given to the
term \lq\lq look like\rq\rq.

It is easy to check that  $\Pb{\pi(v)=v}=1/n$ for all $v\in V$.
Therefore, the expected number of fixed points of $\pi$ is $1$.
However, if $v$ does not appear in any of the
transpositions $T_1,T_2,\dots,T_t$, then $\pi_t(v)=v$.
By the familiar solution of the coupon collector's problem,
we see that when $t=o(n\log n)$, the probability that $\pi_t$ has at most
one fixed point is small. In this sense, $\pi_t$ and $\pi$ are rather
different when $t=o(n\log n)$.
On the other hand, when $t>c\,n\log n$, $c>1/2$, 
the total variation distance between the law of $\pi_t$ and
that of $\pi$ tends to zero as $n\to\infty$~\cite{\DiaconisShahshahani}.

We now consider the situation where $t\le c\,n$ with $c<1/2$.
Let \hfff{d.Gt}{$G^t$} be the graph on $V=\{1,2,\dots,n\}$ where $\{v,u\}$ is an
edge in $G^t$ if and only if the transposition
$(v,u)$ appears in $\{T_1,\dots,T_t\}$.
Let \hfff{d.VGt}{$V_G^t$} denote the set of vertices of the largest connected component of $G^t$
(with arbitrary tie breaking if there is more than one).
By the Erd\H os-R\'enyi Theorem, when $t\le c\,n$, $c<1/2$,
we have $|V_G^t|=O(\log n)$ asymptotically almost surely (a.a.s.).
It follows that the largest cycle (orbit) of $\pi_t$ is also of size $O(\log n)$.
When $c=1/2$, the same holds, but with $\log n$ replaced by any function growing faster
than $n^{2/3}$.
This contrasts with the fact that for every $k\in\{1,2,\dots,n\}$ the
probability that the cycle of $\pi$ containing $1$ has size $\le k$ is precisely $k/n$.
Thus, $\pi_t$ is very different from $\pi$ when $t/n\le c\le 1/2$.

Our main theorem deals with the case where $t/n\ge c>1/2$.
Confirming a conjecture by Aldous, we prove that in this range,
the large cycles of $\pi_t$, when normalized by their total length,
have a distribution that is close to that of
the large cycles of $\pi$. 
A more precise statement of this result will be given shortly.

\medskip

Let $\sigma$ be some permutation of $V$.  
\hfff{d.X}{Let $X(\sigma)$ denote the set of cycles
(orbits) of elements of $V$ under $\sigma$.}  The {\bf cycle structure}
\hfff{d.SS}{
$
\SS(\sigma) 
$}
is then the sorted list of the lengths of the cycles, that is, the list
$\bigl(|C|:C\in X(\sigma)\bigr)$ sorted in nonincreasing order.
Thus,  $\SS_i(\sigma)$ denotes the size of the $i$-th largest cycle of $\sigma$.
If $i$ is larger than the number of cycles of $\sigma$, then we set
$\SS_i(\sigma)=0$, by convention.
Since each $T_j$ is chosen uniformly among the transpositions, it follows that
for each fixed permutation $\sigma$ of $V$ the distribution of $\pi_t$ is the
same as that of $\sigma\circ\pi_t\circ\sigma^{-1}$. Thus, the distribution of $\pi_t$ is
determined by the distribution of the conjugacy class of $\pi_t$.
Now, the conjugacy class of $\pi_t$ is determined by $\SS(\pi_t)$.
Consequently, the distribution of $\SS(\pi_t)$ determines the distribution of
$\pi_t$.

We are now ready to state our main theorem, which gives a positive answer
to a conjecture by David Aldous as stated in~\cite{\BerestyckiDurrett}.

\begin{theorem}\label{t.main}
Let $c>1/2$, and take $t\ge c\,n$.
As $n\to\infty$, the law of $\SS(\pi_t)/|V_G^t|$ converges weakly to
the Poisson-Dirichlet distribution $PD(1)$ with parameter
$1$ (which is defined below).
\end{theorem}

A more explicit statement of the theorem is as follows.
Given $c>1/2$ and $\eps>0$, there is an $n(c,\eps)$
such that for every $n>n(c,\eps)$ and every $t\ge c\,n$ 
there is a coupling of the sequence of transpositions $T_j$ and a
$PD(1)$ sample $Y$ such that
$$
\PB{
\bigl\|Y-\SS(\pi_t)/|V_G^t|\bigr\|_\infty<\eps}>1-\eps\,.
$$

Weak convergence has several equivalent formulations
(see~\cite{\Dudley}), and we have opted to use the coupling version here.

\medskip

The \hfff{d.PD}{$PD(1)$} distribution is a probability
measure on the infinite dimensional simplex
\hfff{d.Om}{
$$
\Omega:=\bigl\{y\in [0,1]^\Np:\sum_iy_i=1,\, y_1\ge y_2\ge \cdots\}
$$}
and may be defined as follows.
Let $U_1,U_2,\dots$ be an i.i.d.\ sequence of random variables uniformly distributed
in $[0,1]$. Set $x_1:=U_1$ and inductively
$x_j:= U_j\bigl(1-\sum_{i=1}^{j-1}x_i\bigr)$.
Let $(y_i)$ be the sequence $(x_i)$ sorted in nonincreasing order.
The $PD(1)$ distribution is defined as the law of $(y_i)$.
See, for example,~\cite{\Holst} for other definitions and a discussion of some
of the properties of the Poisson-Dirichlet distributions.

The behaviour of the size of the largest cluster of $G^t$, which is
the normalizing quantity $|V^t_G|$ in the theorem is known precisely.
The Erd\H os-R\'enyi theorem (see, e.g.,~\cite{\AlonSpencer}) tells us that 
\begin{equation}\label{e.z}
|V_G^t|/n\to z(2\,t/n)
\end{equation}
in probability as $n\to\infty$, where
\hfff{d.z}{$ z(s) $} is the survival probability of a Galton-Watson branching process
with offspring distribution which is the Poisson random variable with mean $s$.
Moreover, $z(s)$ is the positive solution of the equation $1-z=\exp(-s\, z)$
when $s>1$ and $z(s)=0$ for $s\in[0,1]$.

\bigskip

Berestycki and Durrett~\cite{\BerestyckiDurrett} have analysed other aspects
of the chain $\pi_t$ which exhibit a phase transition near $t=n/2$:
they investigate the minimal number of transpositions necessary to write $\pi_t$
as a composition.

In~\cite{\riffle}, Diaconis, McGrath and Pitman discuss the Riffle shuffle,
which is another example where the large cycles appear relaxed well before the permutation
is uniformly distributed.

\bigskip

The evolution of $\SS(\pi_t)$ is also known as the discrete
uniform
coagulation-fragmentation process.
Let us briefly describe the transition from $\SS(\pi_t)$ to $\SS(\pi_{t+1})$.
Suppose that $T_{t+1}$ is the transposition $(a,b)$.
Then $a$ and $b$ are selected uniformly from $V$, and are \lq\lq almost indepedent\rq\rq.
(We could also allow $a=b$, then $T=(a,a)$ would be the identity transposition,
and $a$ and $b$ would be independent. That would not change anything
significant in the following.)
Let $X_i,X_j\in X(\pi_t)$ satisfy $a\in X_i$, $b\in X_j$.
Then $X_i$ and $X_j$ are size biased selections from $X(\pi_t)$,
and are nearly independent given $\pi_t$.
If $X_i\ne X_j$, then in $\pi_{t+1}$ the two cycles
$X_i$ and $X_j$ are replaced by the single cycle whose vertices
are $X_i\cup X_j$. If $X_i=X_j$, then this cycle splits into
two cycles of $\pi_{t+1}$. If $k=|X_i|$ and $m\in\Np$ is the least positive
integer satisfying $\pi_t^m(a)=b$, then the resulting two cycles
of $\pi_{t+1}$ are $\bigl(a,\pi_t(a),\dots,\pi_t^{m-1}(a)\bigr)$
and $\bigl(b,\pi_t(b),\dots,\pi_t^{k-m-1}(b)\bigr)$.
Note that given $X_i$ and given $X_i=X_j$, the resulting
two new cycles have sizes $m$ and $|X_i|-m$, where
$m$ is chosen uniformly in $\{1,2,\dots,|X_i|-1\}$.

There is a similar continuous coagulation-fragmentation
process, 
which is a discrete time Markov chain on the infinite dimensional simplex $\Omega$.
The transition kernel \hfff{d.M}{$M$} of the chain operates as follows.
Given $Y=(Y_1,Y_2,\dots)\in\Omega$,
we choose two indices $i,j\in\Np$ independently, with
$\Pb{i=k\md Y}=\Pb{j=k\md Y}=Y_k$.
If $i\ne j$, then let $Y'$ be obtained from $Y$ by
replacing the two entries $Y_i$ and $Y_j$ with the single
entry $Y_i+Y_j$ and resorting.
If $i=j$, then given $(Y,i,j)$, a random variable 
$v$ is selected uniformly in $[0,Y_i]$ and $Y'$ is obtained
by splitting the entry $Y_i$ into the two entries
$v$ and $Y_i-v$, and resorting. Then $Y'$ is the new state of the
Markov chain.

It is known that the probability measure $PD(1)$ is invariant under $M$.
Apperently, this was first proved in~\cite{\Watterson};
references for several other proofs of this fact are given in~\cite{\DMZZ}.
Vershik conjectured that $PD(1)$ is the only invariant measure.
Subsequently, this was proved by
Diaconis, Mayer-Wolf, Zeitouni and Zerner:

\begin{theorem}[\cite{\DMZZ}]\label{t.DMZZ}
The invariant measure for $M$ is unique.
\end{theorem}

See~\cite{\DMZZ} for more information and
bibliography regarding the history of the problem,
including some earlier established special cases.

The proof of~\cite{\DMZZ} relies on coupling the discrete and the continuous
coagulation-fragmentation processes,
and using representation theory on the symmetric group to understand the
discrete process. In the present paper, we use a different coupling to handle
the continuous process directly, and thereby give a different proof of
Theorem~\ref{t.DMZZ}. Moreover, a slight modification of this coupling
will be essential in the proof of Theorem~\ref{t.main}.

\bigskip

The problems addressed in this paper are a mean-field 
version of a statistical physics model suggested by T\'oth~\cite{\Toth},
which may be described as follows. 
Consider a locally finite graph $G=(V,E)$, and fix a parameter
$\beta>0$. For each (unoriented) edge $e\in E$, let $Z_e\subset[0,1]$ be an independent
Poisson point process of intensity $\beta$ on $[0,1]$.
Let $v_0\in V$. We now describe a walk $v(t)$ starting at
$v(0)=v_0$. Let $t_1$ be the first $t>0$ such that there is
an edge $e_1=[v_0,v_1]$ incident with $v_0$
such that $t_1\in Z_{e_1}+\Z$. If there is no such $t_1$, then
$v(t)=v_0$ for all $t\ge 0$.
But if $t_1$ exists, then let $v(t)=v_0$ for $t\in[0,t_1)$
and $v(t_1)=v_1$.
Inductively, assume that $t_j$ and $v_j$ are defined and
$v(t_j)=v_j$. Let $t_{j+1}$ be the first $t>t_j$
such that there is an edge $e_{j+1}=[v_j,v_{j+1}]$
incident with $v_j$ such that
$t\in Z_{e_{j+1}}+\Z$; set $v(t)=v_j$ for $t\in(t_j,t_{j+1})$
and $v(t_{j+1})=v_{j+1}$.

In the case where $G$ is the complete graph on $V$,
it is easy to see that the orbit of $1$ in $\pi_t$ is
analogous to the range of this walk starting at $1$,
where $\beta=t/n$. The essential difference between
the two is the distinction between continuous time and discrete time.

There are several known open problems regarding T\'oth's model.
Is it true that for (connected) bounded degree graphs $G$, the simple random walk
on $G$ is transient iff T\'oth's walk $v_j$ visits infinitely many vertices
with positive probability for some $\beta>0$?
In particular, is this true for $G=\Z^d$?
For finite graphs $G$, one may ask about the distribution of the size
of the image of the walk $\{v_j\}$, for example.
See~\cite{\Angel} for an analysis of T\'oth's model on regular trees
and for a list of some open problems, including those mentioned above.

Returning to the symmetric group, one may ask about the typical cycle
structure near the transition point $t=n/2$.
A very thorough analogous theory exists for the Erd\H os-R\'enyi transition.
See, for example, \cite{\SpencerLectures,\AlonSpencer,\JLR} and the references cited
there.

\filbreak
\subsection*{Notations}

For the convenience of the reader, we list here some of the notations used extensively,
with hyperlinks and page numbers of the definitions,
and a brief description, where appropriate.

\bigskip
\def\qq{&}

\begin{center}
\halign{
#\quad\hfill&#\quad\hfill&\quad\hfill#\cr
${V}$ \qq $\{1,2,\dots,n\}$ \fff{d.V}
$T_1,T_2,\dots$ \qq i.i.d.\ uniform transpositions on $V$ \fff{d.T}
$\pi_t$ \qq $T_t\circ T_{t-1}\circ\cdots\circ T_1$ \fff{d.pi}
{$X(\sigma)$} \qq set of cycles of a permutation $\sigma$ \fff{d.X}
$X^s$ \qq $X(\pi_s)$ \fff{d.Xt}
{$\SS(\sigma)$} \qq cycle structure of $\sigma$\fff{d.SS}
$X^s(v)$ \qq cycle of $\pi_s$ containing $v$\fff{d.Xsv}
{$V_X^s(k)$} \qq set of cycles of $\pi_s$ of size at least $k$  \fff{d.VX}
$G^t$ \qq graph whose edges correspond to transpositions $T_i,\,i\le t$ \fff{d.Gt}
$V_G^t$ \qq largest cluster in $G^t$ \fff{d.VGt}
$V_G^t(k)$ \qq union of clusters of $G^t$ of size at least $k$ \fff{d.VGs}
$z(s)$ \qq function in the Erd\H os R\'enyi theorem \fff{d.z}
$\Omega$  \qq $\bigl\{y\in [0,1]^\Np:\sum_iy_i=1,\, y_1\ge y_2\ge \cdots\}$\fff{d.Om}
$PD(1)$ \qq Poisson-Dirichlet distribution with parameter $1$\fff{d.PD}
$M$ \qq coagulation-fragmentation transition kernel\fff{d.M}
{$\tilde M$} \qq the coupling \fff{d.tM}
$I(Y,Z)$ \qq indexes of matched entries \fff{d.I}
$Q$ \qq sum of matched entries \fff{d.Q}
$\tilde Y, \tilde Z,\hat Y,\hat Z$ \qq partitions used in defining $\tilde M$ \fff{d.tY}
$u,v$ \qq random variables used in the definition of $\tilde M$ \fff{d.uv}
{$\bar\eps$} \qq $\eps + \text{fragments smaller than $\eps$}$  \fff{d.bareps}
$N^t$ \qq unmatched entries larger than $\eps$ \fff{d.Nt}
$y_1^t, z_1^t$ \qq largest unmatched entries in $Y^t$ and $Z^t$ \fff{d.y1}
}\end{center}

\filbreak

\section{Big pieces}

The main goal of the present section is to show in a quantitative way that most
vertices in $V_G^t$ are in reasonably large cycles of $\pi_t$.

Suppose that $\pi$ is a permutation on $V$ and $T=(x,y)$ a transposition.
If $x$ and $y$ are in different cycles in $\pi$, then
in $T\circ \pi$ these two cycles are joined, and the other cycles remain
unchanged.
Now suppose that $C=(x_0,x_1,\dots,x_m)$ is a cycle of $\pi$
which contains $x$ and $y$. Say, $x=x_j$, $y=x_i$, and $j<i$.
Then in $T\circ\pi$ the cycle $C$ is split into the cycles
$(x_i,x_{i+1},\dots,x_{j-1})$ and
$(x_j,x_{j+1},\dots,x_m,x_0,x_1,\dots,x_{i-1})$.
The other cycles remain unchanged, of course.
This clearly implies the following

\begin{lemma}\label{l.cut}
Let $\pi$ be a permutation of $V$ and $s\in\N$.
Let $T$ be a uniform-random transposition on $V$.
Then the probability that some cycle of $\pi$ is split in
$T\circ\pi$ into two cycles at least one of which has length $\le s$
is at most $2\,s/(n-1)$. \QED
\end{lemma}

This will be used in the next lemma.
Let \hfff{d.Xt}{$X^s=X(\pi_s)$} be the set of cycles of $\pi_s$
and for $v\in V$ let \hfff{d.Xsv}{$X^s(v)$} be the cycle in
$X^s$ containing $v$.
Let \hfff{d.VGs}{$V_G^s(k)\subset V$} be the union of those
connected components of $G^s$ which have at least $k$ vertices,
and let \hfff{d.VX}{$V_X^s(k)\subset V$
be the union of the 
cycles in $X^s$ that have at least $k$ vertices}.

\begin{lemma}\label{l.nofrag}
$$
\E{\bigl|V_G^s(k)\setminus V_X^s(k)\bigr|}\le 4\,s\,k^2/(n-1)
$$
holds for every $k,s\in\N$.
\end{lemma}

\proof
Let $I$ be the set of $t\in\N$ such
that there is a cycle $A\in X^{t-1}$ which splits into two nonempty cycles
in $X^{t}$, $A=A_1\cup A_2$, $A_1,A_2\in X^t$ and at least one of these
cycles, say $A_1$, satisfies $|A_1|\le k$.
The above lemma shows that
$\P[t\in I]\le 2\,k/(n-1)$
for every $t\in\N$,
and hence $\Eb{\bigl|I\cap[0,s]\bigr|}\le 2\,s\, k/(n-1)$.

Suppose that $C\in X^s$, $|C|<k$ and $C\subset V_G^s(k)$.
There must be some vertex $u\in C$ and some
time $t\le s$ such that $|X^t(u)|<|X^{t-1}(u)|$;
otherwise, $C$ would be equal to a component of $G^s$.
Among all such possible pairs $(u,t)$, we choose one that maximizes $t$.
Then we have $X^t(u)\subset C$. Consequently,
$t\in I\cap[0,s]$ and at least one of the two elements of $V$ transposed by $T_t$
is in $C$.
Therefore, the number of such $C$
is at most $2\bigl|I\cap[0,s]\bigr|$.
The statement of the lemma now follows from the above bound
on $\E\bigl|I\cap[0,s]\bigr|$.
\QED

The following lemma will tell us that if $X^t$ has many vertices
in reasonably large cycles at time $t=t_0$, then
with high probability at a specified later time $t_1$
most of these vertices will be in cycles of size at least $\eps\,n$.

\begin{lemma}\label{l.grow}
Let $\delta\in(0,1]$, $t_0,j\in\N$, and $\eps\in(0,1/8)$.
(The lemma will be useful primarily when
$(\log n)^2\le 2^j\le n^{\alpha}$ with any constant $\alpha<1/2$.)
Assume that $2^j<\eps\,\delta\,n$
and that $\Pb{|V_X^{t_0}(2^j)|>\delta\,n}>0$.
Set  $\rho:=2^j/n$ and
\begin{equation}\label{e.twot}
t_1 := t_0+ \bceil{2^6\,\delta^{-1}\, \rho^{-1}\log_2 (\rho^{-1})}\,.
\end{equation}
Then the number of vertices $v$ that are in cycles of size
at least $2^j$ at time $t_0$ but are not in cycles of size
at least $\eps\,\delta\,n$ at time $t_1$ satisfies
\begin{equation}\label{e.failgrow}
\EB{\left|V_X^{t_0}(2^j)\setminus V_X^{t_1}(\eps\,\delta\,n)\right|\md |V_X^{t_0}(2^j)|>\delta\,n}
\le
O(1)\,\delta^{-1}\,\eps\,|\log(\eps\,\delta)|\,n
\,,
\end{equation}
where the constant implied in the $O(1)$ notation is universal.
\end{lemma}

Two important aspects of this lemma are that the right hand
side of~\eref{e.failgrow} does not depend on $j$
and that $t_1$ does not depend on $\eps$.
(However, $t_1-t_0$ depends primarily on $j$ and the right hand
side of~\eref{e.failgrow} depends primarily on $\eps$.)

Before we begin with the actual proof, here is an informal outline.
Let $v\in V_X^{t_0}(2^j)$.
Set $K:=\ceil{ \log_2(\eps\,\delta\,n)}$.
We will choose a sequence of times
$\tau_j,\tau_{j+1},\dots,\tau_K$.
For $s=j,j+1,\dots,K$, 
when $t\in[\tau_s,\tau_{s+1})$ we will \lq\lq expect\rq\rq\
the size of $X^t(v)$ to be at least $2^s$. This can fail in
either of two scenarious: it may happen because a transposition
cuts the cycle of $v$, or it may happen because no transposition
merges the cycle of $v$ with a sufficiently large cycle.
The probabilities for each of these unfortunate situations will
be appropriately estimated.
The choice of the time interval $\tau_{s+1}-\tau_s$ is somewhat delicate.
If it is too long, then perhaps too many cycles will be cut, while
if it is too short, then cycles will not have enough time to merge.
It turns out that
$$
a_s:=2^4\,\delta^{-1}\, 2^{-s}(\ceil{\log_2 n}-s)\,(n-1)\,.
$$
is roughly the right choice, as will become clear in the course of the proof.

\proof
Within the proof below, expectations and probabilities will be
conditioned on $|V_X^{t_0}(2^j)|>\delta\,n$.
Let $K:=\ceil{ \log_2(\eps\,\delta\,n)}$,
and let $a_s$ be as above.
For $s>j$ let $m_s:=\ceil{a_s}$ and set
$m_j:= t_1-t_0-\sum_{s=j+1}^{K-1}m_s$.
Set $\tau_s:=t_0+\sum_{i=j}^{s-1} m_i$.
Note that $\tau_K=t_1$ and $a_s\le m_s\le O(a_s)$ for $s=j,j+1,\dots,K-1$.

Let $s\in\{j,j+1,\dots,K-1\}$ and $t\in\{\tau_s+1,\tau_s+2\dots,\tau_{s+1}\}$.
Define $F^t\subset V$ to be the set of vertices
$v\in V$ such that $|X^t(v)|<|X^{t-1}(v)|$ and $|X^t(v)|<2^{s+1}$.
Lemma~\ref{l.cut} shows that $\E|F^t|\le 2^{2s+4}/(n-1)$.
We also set $\tilde F^t:=\bigcup_{\tau=t_0+1}^t F^\tau$.
Then
$$
\E|\tilde F^{t_1}|
\le \sum_{s=j}^{K-1} m_s\,2^{2s+4}/(n-1)
= O(\eps\,|\log(\eps\,\delta)|\,n)\,.
$$

We consider the vertices in $F^t$ as vertices \lq\lq failing\rq\rq\ at
time $t$.
However, there are other ways in which vertices can fail.
If at time $t\in\{\tau_s,\tau_s+1,\dots,\tau_{s+1}-1\}$
we have $|V_X^t(2^s)|<\delta\,n/2$, then we consider the
whole process as failed, and we set $H^t:=V$.
Otherwise, take $H^t=\emptyset$.
Also set $\tilde H^t:=\bigcup_{t'=t_0}^t H^{t'}$.

The third and last way in which a vertex $v$ may fail is
if $X^t(v)$ does not grow in time.
Let 
$$
B^s:= V_X^{\tau_s}(2^s)\setminus\Bigl(\tilde F^{\tau_{s+1}}\cup
\tilde H^{\tau_{s+1}-1}\cup V_X^{\tau_{s+1}}(2^{s+1})\Bigr)\,,
$$
and $\tilde B^s:=\bigcup_{k=j}^s B^k$.
The vertices in $B^s$ are vertices whose cycles failed to
grow sufficiently between time $\tau_s$ and time $\tau_{s+1}$.
It is clear that 
\begin{equation}\label{e.cover}
V_X^{t_0}(2^j)\subset V_X^{t_1}(\eps\,\delta\,n)
\cup
\tilde H^{t_1}\cup\tilde F^{t_1}\cup
\tilde B^K
\,.
\end{equation}

If $v\in B^s$, then it must be the case that for every
$t\in[\tau_s,\tau_{s+1}-1]$ we have $|V_X^{t}(2^s)|\ge \delta\,n/2$
(since $B^s$ is disjoint from $\tilde H^{\tau_{s+1}-1}$)
and $v\in V_X^t(2^s)\setminus V_X^t(2^{s+1})$
(since $B^s$ is disjoint from $\tilde F^{\tau_{s+1}}$).
If we condition on $2^s\le |X^t(v)|<2^{s+1}$
and on $|V_X^{t}(2^s)|>\delta\,n/2$, then there is probability
at least 
$$
2^s(\delta\,n/2-2^{s+1})\,{n\choose 2}^{-1}
\ge 2^{s-3} \delta\,(n-1)^{-1}
$$
that $T_{t+1}$ transposes an element from $X^t(v)$
and an element from some other cycle of $X^t$ whose size is at least $2^s$.
If that happens, then $v\in V_X^t(2^{s+1})$ and this implies that 
$v$ cannot be in $B^s$.
Consequently, 
\begin{equation*}
\begin{aligned}
\Pb{v\in B^s } 
&
\le
\Bigl(1- 2^{s-3} \delta/(n-1)\Bigr) ^{m_{s}}
\cr\\&
\le\exp\bigl(- 2^{s-3} \delta\,{m_{s}}/(n-1)\Bigr) 
\le O(2^{s}/n)\,.
\end{aligned}
\end{equation*}
Hence,
$$
\E|\tilde B^K|\le O(1)\,n\,2^{K}\,n^{-1}= O(\eps\,\delta\,n).
$$

It follows from the definition of $H^{t}$ that in order for
$H^{t}$ to be nonempty, we must have
$|\tilde F^t\cup \tilde B^{s-1}|\ge \delta\,n/2$.
Therefore,
$$
\E|\tilde H^{t_1}|\le n\,\Pb{|\tilde F^{t_1}\cup \tilde B^{K}|\ge\delta\,n/2}
\le
2\,\delta^{-1}\E|\tilde F^{t_1}\cup \tilde B^{K}|\,.
$$
When we combine this with~\eref{e.cover} and the above estimates for
$\E|\tilde F^{t_1}|$ and $\E|\tilde B^{K}|$,
the lemma follows.
\QED

\begin{lemma}\label{l.precoup}
Fix some $c>1/2$, and let $t \ge c\, n$, $t\in\N$. 
Let $\eps,\alpha\in(0,1/8)$ and let $N$ be the minimal number of cycles in $X^t$
which cover at least $(1-\eps)\,|V_G^t|$ vertices of $V_G^t$.
Then 
$$
\Pb{N> \alpha^{-1}\,|\log(\alpha\,\eps)|^2}\le C_1\,\alpha
$$
for all $n>n_1$,
where $C_1$ is a constant which depends only on $c$, and $n_1$ may depend on $c$ and $\eps$.
\end{lemma}

\proof
First, suppose that $t\le n^{5/4}$.
Choose $j$ such that $n^{1/4}\le 2^j < 2\,n^{1/4}$.
Let $\delta=z/2$, where $z=z(2\,t/n)$ is the Galton-Watson survival probability
discussed in the introduction.
Choose $t_0$ so that~\eref{e.twot} holds with $t$ in place of $t_1$.
Note that $t-t_0=O(n^{3/4}\log n)$.
(Here and below, the constants in the $O(\cdot)$ notation may depend on $c$.)
We apply the Erd\H os-R\'enyi theorem at time $t_0$ to conclude
that a.a.s.\ $|V_G^{t_0}|- n \,z=o(n)$ and the second
largest component of $G^{t_0}$
has size less than $(\log n)^2$.
Lemma~\ref{l.nofrag} with $k=2^j$ and $s=t_0$ implies that
$|V_G^{t_0}\setminus V_X^{t_0}(2^j)|\le n^{7/8}$ a.a.s.
Note also that $|V_G^{t}\setminus V_G^{t_0}|\le (t-t_0)\,O(\log n)^2$
a.a.s., because we know that the components of $G^{t_0}$ other than the
largest one are typically smaller than $(\log n)^2$.
Hence $|V_G^{t}\setminus V_X^{t_0}(2^j)|<n^{7/8}$ a.a.s.
Now, Lemma~\ref{l.grow} implies that for every fixed $\eps'>0$ and for every sufficiently
large $n$
\begin{equation}\label{e.capture}
\EB{\bigl|V_G^{t}\setminus V_X^{t}(\eps' n)\bigr|} < O(1) \, \eps'\,| \log\eps'|\,n
\,.
\end{equation}
(Note that $|V_G^{t_0}|\le n$, and hence the conditioning in~\eref{e.failgrow}
may be ignored once $n$ is large enough so that
$\Pb{|V_X^{t_0}(2^j)|\le \delta\,n}< \eps'\,|\log\eps'|$.)

Now, to show that~\eref{e.capture} holds also without the assumption that $t\le n^{5/4}$,
we note that Lemma~\ref{l.grow} may be applied with $j=0$, $\delta=1$ and 
$t_0$ chosen so that~\eref{e.twot} holds with $t$ in place of $t_1$. (In this case, we do not
need to use Lemma~\ref{l.nofrag}.)

Set $a(k):=\bigl|V_G^{t}\setminus V_X^{t}(k)\bigr|$.
Let $i_0$ be the smallest integer $i$ such that $a(2^{-i}n)<\eps\,n/2$.
Then $N$ is bounded by the number of cycles in $V_X^{t}(2^{-i_0}n)\cap V_G^t$.
Let $i_1$ be the least integer such that $2^{-i_1}<\alpha\,\eps/|\log (\alpha\,\eps)|$.
Then~\eref{e.capture} shows that $\P[i_0>i_1]=O(\alpha)$.
We may write 
$$
a(k)=\sum\bigl\{|A|:A\subset V_G^t,\,A\in X^t,\,|A|<k\bigr\}\,.
$$
By considering the contribution of each cycle to the sum
$$
S_m:=\sum_{i=0}^{m} a(2^{-i}n)\, 2^i/n
$$
we find that $N=O(S_{i_0})$.
On the other hand~\eref{e.capture} implies that
$$
\E[S_{i_1}]\le O(1) \sum_{i=0}^{i_1} i\le O(1)\,(i_1)^2\le O(1)\,|\log(\alpha\,\eps)|^2\,.
$$
Because $\P[i_0>i_1]=O(\alpha)$, this completes the proof.
\QED

\section{Coupling}\label{s.coup}

At this point, it seems likely that the proof of Theorem~\ref{t.main} 
can be completed using some of the results from the work of
Diaconis, Mayer-Wolf, Zeitouni and Zerner~\cite{\DMZZ}.
However, we prefer instead to use a different coupling argument to finish
off the proof and also prove the main result of~\cite{\DMZZ}.

We now describe a coupling in the continuous setting. A similar coupling
will also apply to couple between the discrete and continuous
setting, but the purely continuous setting avoids several annoying
minor notational issues.

The coupling is between two Markov chains
$Y^t$ and $Z^t$ starting at
possibly different initial starting points $Y^0,Z^0\in\Omega$
with each separately evolving according to the transition kernel $M$.
In this coupling, the evolution of $(Y^t,Z^t)$ will also be Markov.
\hfff{d.tM}{Its transition kernel will be denoted by $\tilde M$.}

The basic idea in the construction of $\tilde M$ is that if we have
entries in $Y^t$ that are equal to entries in $Z^t$, then
we don't want to ruin this. Consequently, if we make a change to such
an entry in $Y^t$, we want to make a corresponding change to the
corresponding entry in $Z^t$. On the other hand, as much as we can,
we do want to produce new entries in $Y^t$ and $Z^t$ that match.
Our measure of the discrepancy between $Y^t$ and $Z^t$ will roughly
be the number of large unmatched entries, and we will strive to
reduce the discrepancy.

In order to define $\tilde M$, we need some more notations.
Let $(Y,Z)\in\Omega^2$.
We will need to match entries in $Y$ with entries in $Z$ of the
same length, if such exist, and match as many
entries as possible.  The matching will be encoded
via maps $f_{Z,Y},f_{Y,Z}:\Np\to\N$, which are defined as follows.
Let $i\in \Np$, let $H$ be the set of $j\in\Np$ such that $Y_i=Z_j$,
and let $k:=\bigl|\{j\in\N:j\le i,\, Y_j=Y_i\}\bigr|$.
(Partly because we want to easily generalize to the discrete setting,
we do not want to rule out the possibility that $Y_i=Y_j$ for some $i\ne j$.)
If $|H|<k$, then set $f_{Y,Z}(i)=0$.
Otherwise, let $f_{Y,Z}(i)$ be the $k$'th smallest element in $H$.
(By exchanging $Y$ and $Z$, this also defines the map $f_{Z,Y}$.)
Let 
$$
\hfff{d.I}{I(Y,Z):=f_{Y,Z}^{-1}(\Np)=\{i\in\N:f_{Y,Z}(i)\ne0\}\,.}
$$
The entries $Y_i$ with $i\in I(Y,Z)$ will be referred to as matched.
Likewise, $Z_j$, $j\in I(Z,Y)$, are the matched entries of $Z$.
Observe that $f_{Z,Y}\circ f_{Y,Z}(i)=i$ for every $i\in I(Y,Z)$,
$f_{Y,Z}(I(Y,Z))=I(Z,Y)$ and
$Z_{f_{Y,Z}(i)}=Y_i$ for every $i\in I(Y,Z)$.
Let 
$$
\hfff{d.Q}{Q=Q(Y,Z):=\sum\{Y_i:i\in I(Y,Z)\}=\sum\{Z_j:j\in I(Z,Y)\}\,.}
$$

\begin{figure}
\medskip
\SetLabels
\R(-0.01*.8)$\tilde Y$\\
\R(-0.01*.3)$\tilde Z$\\
\T(.47*-.01)$u$\\
\T(.62*-.01)$1-Q$\\
\endSetLabels
\centerline{\epsfxsize=4.4in\AffixLabels{\epsfbox{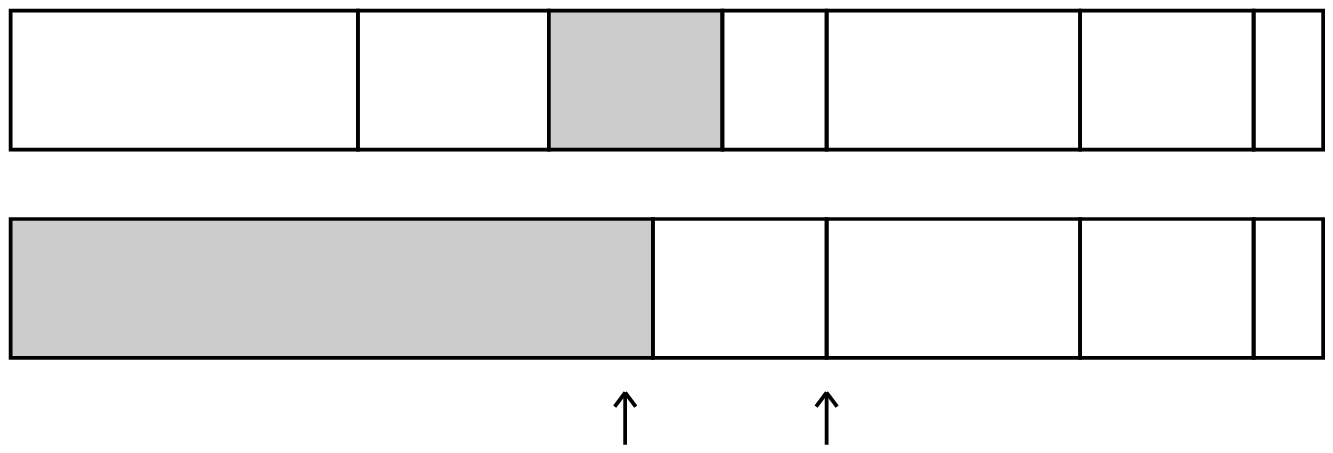}}}
\vskip 0.3in
\SetLabels
\R(-0.01*.8)$\hat Y$\\
\R(-0.01*.3)$\hat Z$\\
\T(.085*-.01)$v$?\\
\T(.43*-.01)$v$?\\
\T(.59*-.01)$v$?\\
\T(.93*-.01)$v$?\\
\endSetLabels
\centerline{\epsfxsize=4.4in\AffixLabels{\epsfbox{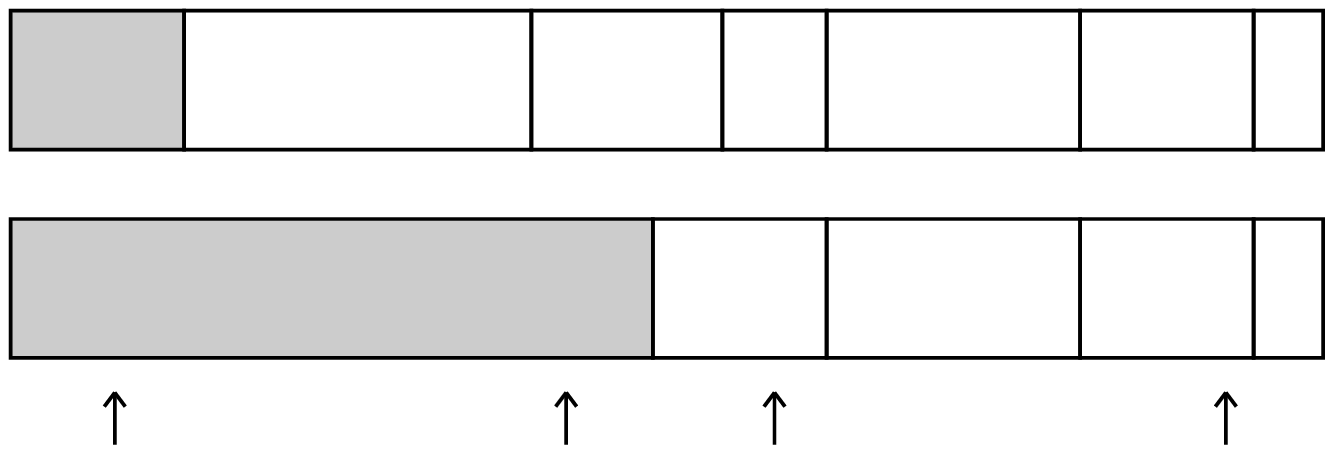}}}
\medskip
\begin{caption}{
\label{f.coup}
The random variable $u$ chooses
a segment in $\tilde Y$ and a segment in $\tilde Z$.
The different illustrated choices for the
random variable $v$ yield a split in $Y$ and in $Z$, a merge in $Y$
and a split in $Z$, a merge in both where a matched segment is not involved,
and merges involving matched segments, respectively.}
\end{caption}
\end{figure}

We will now describe the transition kernel $\tilde M$.
Given $(Y,Z)\in\Omega$,
we need to perform one step of $M$ for each of $Y$ and $Z$,
thereby generating new configurations $Y'$ and $Z'$.
We associate with $Y$ and with $Z$  partitions
\hfff{d.tY}{$\tilde Y=(\tilde Y_i:i\in\Np)$}
and
$\tilde Z=(\tilde Z_i:i\in\Np)$ 
of $[0,1]$ into closed intervals,
as follows. (See also Figure~\ref{f.coup}.)
The length of the interval $\tilde Y_i$
is $Y_i$. The intervals $\tilde Y_i$ with $i\in I(Y,Z)$
tile the interval $[1-Q,1]$, while the intervals with $i\notin I(Y,Z)$
tile the interval $[0,1-Q]$. Within each of these classes, let the
intervals be ordered according to the indices; that is
$\max \tilde Y_i\le \min\tilde Y_{i'}$ if $i<i'$
when $i,i'\in I(Y,Z)$ and when $i,i'\notin I(Y,Z)$.
A partition $\tilde Z=(\tilde Z_j:j\in\Np)$ is
constructed in the same way. Note that necessarily
$\tilde Z_{f_{Y,Z}(i)}=\tilde Y_i$ whenever $i\in I(Y,Z)$.

\hfff{d.uv}{Let $u$ and $v$ be two independent uniform random variables in $[0,1]$.}
Let $a,a'\in\Np$ be the indices satisfying $u\in\tilde Y_a$
and $u\in\tilde Z_{a'}$. In this way, $u$ induces a size biased
sample from $Y$ and from $Z$. We will use $v$ to induce a different
size biased sample, based on different tilings of $[0,1]$.
Let $\hat Y$ be the tiling $(\hat Y_i:i\in\Np)$
of $[0,1]$ by intervals that is obtained from
$\tilde Y$ by shifting the interval $\tilde Y_a$ to the begining.
(That is, $\hat Y_a=[0,Y_a]$,
$\hat Y_i=Y_a+\tilde Y_i$ if $\max \tilde Y_i\le\min\tilde Y_a$
and $\hat Y_i=\tilde Y_i$ if $\max\tilde Y_a\le\min\tilde Y_i$.)
Similarly, $\hat Z$ is the tiling obtained from $\tilde Z$
by shifting $\tilde Z_{a'}$ to the begining.
Note that $\hat Z_{f_{Y,Z}(i)}=\hat Y_i$ whenever $i\in I(Y,Z)$.

Let $b$ and $b'$ be the indices satisfying
$v\in \hat Y_b$ and $v\in\hat Z_{b'}$.
If $a\ne b$, let $Y'$ be obtained from $Y$ be
replacing the two entries $Y_a$ and $Y_{b}$ by the
single entry $Y_a+Y_{b}$ and resorting.
If $a=b$, let $Y'$ be obtained from $Y$ by
replacing $Y_a$ with the two entries
$v$ and $Y_a-v$ and resorting.
Similarly, if $a'\ne b'$, let $Z'$ be obtained from $Z$ be
replacing the two entries $Z_{a'}$ and $Z_{b'}$ by the
single entry $Z_{a'}+Z_{b'}$ and resorting.
If $a'=b'$, let $Z'$ be obtained from $Z$ by
replacing $Z_{a'}$ with the two entries
$v$ and $Z_{a'}-v$ and resorting.
This completes the construction of the Markov transition kernel
$\tilde M$.

Let us observe a few essential features of this coupling.
If $Y_i$ and $Z_j$ are split, then one of the two new entries
in each of $Y'$ and $Z'$ is equal to $v$.
If $i\in I(Y,Z)$ and $Y_i$ is split or merged, then the same
happens to $Z_{f_{Y,Z}(i)}$.
Similarly, if $j\in I(Z,Y)$ and $Z_j$ is split or merged, then
the same happens to $Y_{f_{Z,Y}(j)}$.

\bigskip

We first informally describe the general behaviour of $\tilde M$, postponing the
exact statements and proofs.
When there are several unmatched reasonably large entries in $Y^t$ and in $Z^t$,
these merge and become few quite quickly.
However, when they are very few, it is hard for them to dissappear completely.
Suppose that there is
one large unmatched entry in $Y^t$ and two unmatched entries in $Z^t$.
When the two unmatched entries in $Z^t$ are merged, the single unmatched entry
in $Y^t$ is likely to be split. Thus, the situation does not improve so quickly.
There is a parity phenomenon here: if the number of positive entries in $Y^t$ is finite,
then its parity either stays the same as that of $t$, or is opposite to that of $t$.
Even if the number of positive entries is infinite, if it takes a long time for the
smaller entries to be hit, the larger entries appear to follow this parity periodicity.
One way to handle the parity issue would be to introduce a delay to either $Y^t$ or $Z^t$,
but not both, in order to match up their parities. However,  another phenomenon
will be used instead. An unmatched entry in $Y^t$ often splits into one matched
entry and one unmatched entry. With any luck, the unmatched entry might be rather
small. Thus, large unmatched entries are replaced by small unmatched entries.
In effect, there is a diffusion of unmatched entries between different scales.
Because of this, it is eventually unlikely to find a large unmatched entry,
which is what we want to prove. However, this latter process is much slower than the
first stage where large unmatched entries merge and become fewer. Thus,
in time $t$ the largest unmatched entry one can expect to find is of order roughly
$1/\log t$.

\bigskip

Define
$$
N_\eps(Y,Z):= \bigl|\{i\in \Np\setminus I(Y,Z): Y_i>\eps\}\bigr|
.
$$
This is the number of entries in $Y$ that are not matched by entries
in $Z$ and have size larger than $\eps$.

\begin{lemma}\label{l.coup1}
Let $\eps>0$, and let $Y^0,Z^0\in\Omega$.
Let $(Y^t,Z^t)$ be the Markov chain given by $\tilde M$ starting
at $(Y^0,Z^0)$.
To abbreviate notations, set
\hfff{d.Nt}{$N^t:=N_\eps(Y^t,Z^t)+N_\eps(Z^t,Y^t)$},
$Q^t=Q(Y^t,Z^t)=Q(Z^t,Y^t)$,
$I^t:=I(Y^t,Z^t)$ and $J^t:= I(Z^t,Y^t)$.
Also define
\hfff{d.bareps}{
$$
\bar\eps:=\eps+
\sum \{Y^0_i:Y^0_i<\epsilon\}+
\sum \{Z^0_i:Z^0_i<\epsilon\}\,.
$$}
Let \hfff{d.y1}{$y_1^t:=\max\bigl\{Y^t_i:i\notin I^t\bigr\}$}
be the size of the largest unmatched entry of $Y^t$
(set $y_1^t=0$ if all entries are matched),
and let $z_1^t$ be the size of the largest unmatched entry of $Z^t$.
Let $\TT$ be a random variable with values in $\N$ which is independent
from the evolution of the chain $(Y^t,Z^t)$.
Set 
$$
\eta:=\max\bigl \{ \P[\TT=t]:t\in\N\bigr\}\,.
$$
Then  
\begin{equation}\label{e.yz}
\EB{(1-Q^\TT)\bigl(1-Q^\TT-\max\{y_1^\TT,z_1^\TT\}\bigr)} \le 
\frac\eta2\,N^0+4\,\bar\eps\,\E[\TT+1]\,.
\end{equation}
\end{lemma}

When the right hand side in~\eref{e.yz} is small, we know that
with high probability either the sum
of the unmatched entries in $Y^\TT$ is only slightly larger than the
largest unmatched entry, or this is true for $Z^\TT$.

\proof
Let $\ev A_s$ be the event
that up to time $s$ in every merging occuring both merged pieces are of size
at least $\eps$ and in every splitting both resulting pieces are of size
at least $\eps$.
Let $\ev F_s$ be the $\sigma$-field generated by $\bigl((Y^t,Z^t):t=0,1,\dots,s\bigr)$.
Conditioned on $\ev F_{t-1}$, the
probability that at time $t$ there is a split in any $Y^{t-1}_i$
and one of the pieces is of size less than $\eps$ is at most $2\,\eps$.
Conditioned on $\ev F_{t-1}$, the
probability that there is any $Y^{t-1}_i$ with
$Y^{t-1}_i<\eps$ that is merged at time $t$ with some other $Y^{t-1}_{j}$ is at most
$2\sum\{Y_i^{t-1}:Y_i^{t-1}<\eps\}$.
Similar considerations apply to $Z^t$.  Consequently, for $t\in\Np$,
\begin{equation}
\label{e.A}
\Pb{\neg\ev A_{t}\md \ev A_{t-1},\,\ev F_{t-1}}\le 4\,\bar\eps\,.
\end{equation}

We now study the evolution of the quantity $N^t$, and consider
several different cases for the transition from
$(Y^{t},Z^{t})$ to $(Y^{t+1},Z^{t+1})$. In each case we assume that $\ev A_{t+1}$ holds.
\begin{enumerate}
\item The transition involves splitting in $Y^{t}$ and merging in $Z^{t}$.
Suppose that $Y_i^{t}$ is split and $Z^{t}_j$ is merged with
$Z^{t}_{j'}$. Then necessarily $i\notin I^t$ and $j,j'\notin J^t$.
Since $\ev A_{t+1}$ is assumed to hold, it follows that
$N_\eps(Y^{t+1},Z^{t+1})\le N_\eps(Y^t,Z^t)+1$ and
$N_\eps(Z^{t+1},Y^{t+1})\le N_\eps(Z^t,Y^t)-1$.
Thus, in this case, $N^{t+1}\le N^t$.
\item The transition involves splitting in $Z^t$ and merging in $Y^t$.
By symmetry, also in this case we have $N^{t+1}\le N^t$.
\item The transition involves splitting in $Y^t$ and splitting in $Z^t$.
Note that by construction the size of one of the newly created split
entries is the same for $Y$ as for $Z$.
Suppose that $Y_i^t$  and $Z^t_j$ are split. 
If $i\in I^t$ then also $j\in J^t$ and $Y^t_i=Z^t_j$.
In that case, both new entries for $Y$ are the same as the new entries
for $Z$, and hence $N^{t+1}=N^t$.  The same conclusion is obtained if
$j\in J^t$. If $i\notin I^t$ and $j\notin J^t$,
then in both $Y^t$ and $Z^t$ an unmatched entry is replaced
by two entries at least one of which is matched. Thus $N^{t+1}\le N^t$.
\item The transition involves merging in $Z^t$ and merging in $Y^t$.
Suppose that $Y^t_i$ is merged with $Y^t_{i'}$.
It is easy to verify, as above, that in this case also $N^{t+1}\le N^t$.
However, if $i,i'\notin I^t$, then the corresponding statement
is also true for the merged entries in $Z^t$, and we actually have
$N^{t+1}\le N^t-2$.
\end{enumerate}
In summary, we see that on the event $\ev A_{t+1}$ we
have $N^{t+1}\le N^t$ and $N^{t+1}\le N^t-2$ when
there is merging in both $Y^t$ and $Z^t$ and the merging
does not involve matched entries.

Since $N^t\ge 0$, we obviously have
$$
\sum_{t=0}^\infty (N^t-N^{t+1})\,1_{\ev A_{t+1}}\le N^0\,,
$$
and we have seen that all the summands are nonnegative.
Since $\TT$ is independent from $(N^t-N^{t+1})1_{\ev A_{t+1}}$,
\begin{equation}\label{e.unic}
\begin{aligned}
&
\Eb{(N^\TT-N^{\TT+1})\, 1_{\ev A_{\TT+1}}}
= \sum_t \Eb{ (N^t-N^{t+1})\, 1_{\ev A_{t+1}}\,1_{\TT=t}}
\cr & \qquad\qquad\qquad\qquad
= \sum_t \Eb{ (N^t-N^{t+1})\, 1_{\ev A_{t+1}}}\,\P[\TT=t]
\le \eta\,N^0\,.
\end{aligned}
\end{equation}

Set $a^t=1-Q^t-\max\{y_1^t,z_1^t\}$.
Recall the random variables $u$ and $v$ used in the transition
kernel $\tilde M$.
If in the transition from $(Y^t,Z^t)$ to $(Y^{t+1},Z^{t+1})$
we have 
$u < 1-Q^t$ and $\max\{y_1^t,z_1^t\}< v< 1-Q^t$,
then in both $Y$ and $Z$ we have merging
of unmatched entries.
Thus, 
$$
\PB{N^t-N^{t+1}\ge 2\text{ or }{\neg \ev A_{t+1}}\md
\ev F_t} \ge
(1-Q^t)\,a^t\,.
$$
By applying this at time $t=\TT$ and taking expectations, we get
\begin{equation*}
\begin{aligned}
\Eb{ (1-Q^\TT)\,a^\TT }
&
\le
\Pb{N^\TT-N^{\TT+1}\ge 2\text{ or }{\neg \ev A_{\TT+1}}}
\cr &
\le
\frac 12\,\Eb{(N^\TT-N^{\TT+1})\,1_{\ev A_{\TT+1}}}+ \Pb{\neg \ev A_{\TT+1}}
\,.
\end{aligned}
\end{equation*}
Consequently,~\eref{e.A} and~\eref{e.unic} complete the proof
of the lemma.
\QED

Assuming that we can make the right hand side of~\eref{e.yz}
small, Lemma~\ref{l.coup1} tells us that with high probability
either $1-Q^\TT-y^\TT_1$ or
$1-Q^\TT-z^\TT_1$ is small. 
If we knew that both are small, it would
follow that also $y_1^\TT-z_1^\TT$ is rather small, since 
$\sum_i Y_i^\TT=\sum_j Z_j^\TT=1$.
However, it might be the case that $1-Q^\TT-z^\TT_1$ is small
but $1-Q^\TT-y^\TT_1$ is not.  The next lemma tells us that in such a situation,
with high probability, $Y^\TT$ does not have more than two significant
unmatched entries.

\begin{lemma}\label{l.coup2}
With the setting and notations of Lemma~\ref{l.coup1},
let $y^t_2$ be the second largest unmatched entry in $Y^t$.
For every $\rho\in(0,1)$
\begin{equation}\label{e.coup2}
\Pb{1-Q^\TT-y^\TT_1-y^\TT_2>\rho}< 
2^6\,\rho^{-4}\, \eta\,N^0+2^9\,\bar\eps\,\rho^{-4}\,\E[\TT+2].
\end{equation}
\end{lemma}
\proof
Let $\ev D$ be the event $\{1-Q^\TT-y^\TT_1-y^\TT_2>\rho\}$ and let
$\ev R$ be the event $\{1-Q^\TT-z^\TT_1<\rho/4\}$.
Assume that $\ev D\cap\ev R$ holds.
Then $z^\TT_1\ge 3\rho/4+y^\TT_1+y^\TT_2$.
Let $\ev U$ be the event that the random variables $u$ and $v$
used in the transition from 
$(Y^\TT,Z^\TT)$ to $(Y^{\TT+1},Z^{\TT+1})$
satisfy $u<3\rho/4$ and $z^\TT_1-\rho/2<v<z^\TT_1-\rho/4$.
On $\ev D\cap\ev R\cap\ev U$,
the largest unmatched entry in $Z^\TT$ will be split
and the transition from $Y^\TT$ to $Y^{\TT+1}$ would involve
a merge (of unmatched entries), because $z^\TT_1-\rho/2>y_1$.
Consequently, a.s.\ on $\ev D\cap\ev R\cap\ev U$
 the two new entries of
$Z^{\TT+1}$ will be unmatched in $Y^{\TT+1}$, and
in particular, $1-Q^{\TT+1}\ge z^\TT_1$.
Moreover, each of the new entries of $Z^{\TT+1}$ would be larger than
$\rho/4$. Clearly, $y^{\TT+1}_1\le y^\TT_1+y^\TT_2$.
Consequently, $1-Q^{\TT+1}-y_1^{\TT+1}\ge z^\TT_1-y^\TT_1-y^\TT_2\ge 3\rho/4$
and $1-Q^{\TT+1}-z^{\TT+1}_1\ge \rho/4$.
Thus, on $\ev D\cap\ev R\cap\ev U$, we have
$1-Q^{\TT+1}-\max\{y_1^{\TT+1},z_1^{\TT+1}\}\ge\rho/4$.
Now,
\begin{equation*}\begin{aligned}
&
\EB{(1-Q^{\TT+1})\bigl(1-Q^{\TT+1}-\max\{y_1^{\TT+1},z_1^{\TT+1}\}\bigr)} 
\cr &\qquad
\ge 
\EB{\bigl(1-Q^{\TT+1}-\max\{y_1^{\TT+1},z_1^{\TT+1}\}\bigr)^2\md \ev D,\,\ev R,\,\ev U}\,
\Pb{\ev D,\,\ev R,\,\ev U}
\cr &\qquad
\ge 
(\rho^2/16)\,
\Pb{\ev D,\,\ev R,\,\ev U}
.
\end{aligned}\end{equation*}
Lemma~\ref{l.coup1} with $\TT$ replaced by $\TT+1$ therefore gives
$$
\Pb{\ev D,\,\ev R,\,\ev U}
\le
8\,\rho^{-2}
\eta\,N^0+2^6\,\bar\eps\,\rho^{-2}\,\E[\TT+2].
$$
Clearly, $\Pb{\ev U\md\ev D,\,\ev R}=3\,\rho^2/16$, and hence
$$
\Pb{\ev D,\,\ev R}\le (16/3)\,\rho^{-2}\,\Pb{\ev D,\,\ev R,\,\ev U}.
$$
On the other hand, on $\ev D\setminus\ev R$ we have
$(1-Q^\TT)\bigl(1-Q^\TT-\max\{y_1^\TT,z_1^\TT\}\bigr)> \rho^2/4$.
Thus, applying Lemma~\ref{l.coup1} again gives
$$
\Pb{\ev D\setminus\ev R} \le
2\,\rho^{-2}\,\eta\,N^0+16\,\bar\eps\,\rho^{-2}\,\E[\TT+1].
$$
Since $ \Ps{\ev D}=\Ps{\ev D,\,\ev R}+\Ps{\ev D\setminus\ev R}$,
the above estimates combine to give~\eref{e.coup2}, and complete the proof.
\QED


\begin{lemma}\label{l.coup3}
With the setting and notations of Lemma~\ref{l.coup1},
Let $\rho\in(0,1/8)$ and assume that $0<\eps<\rho$.
Then for each $t\in\Np$ and for every $n\in\Np$ satisfying $2^n\le t\,\rho$
\begin{equation}
\label{e.coup3}
t^{-1}\sum_{\tau=0}^{t-1} \Pb{y_1^\tau\ge \rho }
\le
 O(\rho^{-1}n^{-1})+
O(2^{4n}/\rho^{5})(N^0/t+\bar\eps\,t).
\end{equation}
\end{lemma}

The basic idea of the proof of the lemma is to use the fact that conditioned
on $y_1^\tau\ge \rho$ there is a significant enough probability that at a later
time $\sigma$ there will be some unmatched
$Y^\sigma_{i'}\in [2^{-k}\rho,2^{-k+1}\rho]$, since the unmatched
piece at time $\tau$ of size $\ge\rho$ may be split immediately.
Lemma~\ref{l.coup2} is then used to show that when we fix $\sigma$, 
with high probability the latter event occurs for at most three different
$k$ in the range $\{1,\dots,n\}$, if $n$ is not too large. An appropriate
summation over $k$ and $\sigma$ completes the proof.

\proof
For $\sigma>\tau$, $\sigma,\tau\in\N$, $k\in\Np$,
let $\ev X(\tau,\sigma,k)$ be the event that the transition between
time $\tau$ and $\tau+1$ produces a splitting in $Y^\tau$ and one of the
split pieces is unmatched, has size in the range $[2^{-k-1} \rho,2^{-k}\rho)$,
and this split piece is not modified up to time $\sigma$.
Set 
$$
\ev X'(\tau,\sigma,k):=
\ev X(\tau,\sigma,k)\setminus\bigcup_{\tau'=\tau+1}^{\sigma-1}
\ev X(\tau',\sigma,k)\,.
$$
Suppose that $y_1^\tau\ge\rho $ and that $Y_i^\tau=y_1^\tau$.
 If in the transition from
$\tau$ to $\tau+1$ we have $u\in \tilde Y_i^\tau$ 
and $v\in (Y^\tau_i-2^{-k}\rho,Y^\tau_i-2^{-k-1}\rho)$,
then $Y^\tau_i$ is indeed split, and it is easy to see that
the resulting piece of $Y^\tau_i-v$ is unmatched a.s.
If that happens, the conditioned probability that
up to time $\sigma$ this piece is modified is bounded by
$2\,(\sigma-\tau)\,2^{-k}\rho$, since the size of this piece
is at most $2^{-k}\rho$. 
This gives
$$
\Pb{\ev X(\tau,\sigma,k)\md y_1^\tau\ge \rho }
\ge 2^{-k-1}\rho^2\,\bigl(1-2\,(\sigma-\tau)\,2^{-k}\rho\bigr).
$$
On the other hand, the conditional probability for
$\ev X(\tau',\sigma,k)$ given the configuration at time $\tau'$
is clearly at most $2^{-k}\rho$.
Hence
\begin{equation*}
\begin{aligned}
\Pb{\ev X'(\tau,\sigma,k)\md y_1^\tau\ge \rho }
&
\ge
\Pb{\ev X(\tau,\sigma,k)\md y_1^\tau\ge \rho }
\bigl(1-(\sigma-\tau)\,2^{-k}\rho\bigr)
\\ &
\ge
2^{-k-1}\rho^2\,\bigl(1-2\,(\sigma-\tau)\,2^{-k}\rho\bigr)
\bigl(1-(\sigma-\tau)\,2^{-k}\rho\bigr)
\\ &
\ge
2^{-k-1}\rho^2\,\bigl(1-3\,(\sigma-\tau)\,2^{-k}\rho\bigr),
\end{aligned}
\end{equation*}
which implies 
\begin{equation}
\Pb{y_1^\tau\ge\rho}\le
 2^{k+3}\rho^{-2}\,
\Ps{\ev X'(\tau,\sigma,k)},
\qquad\text{if }
\tau<\sigma\le\tau+ 2^{k-2}/\rho\,.
\label{e.evX}
\end{equation}
Let $\ev V^\sigma(n)$ be the event that there are at least $3$ distinct
$k\in\{0,1,\dots,n-1\}$
such that there is an unmatched $Y_i^\sigma$ in the
range $[2^{-k-1}\rho,2^{-k}\rho)$.
We now apply Lemma~\ref{l.coup2} with $\TT$ chosen uniformly in
$\{0,1,\dots,2t-1\}$ and with $\rho$ replaced by $2^{-n}\rho$ to get
\begin{equation}
\label{e.fsum}
\sum_{\sigma=0}^{2t-1} \Ps{\ev V^\sigma(n)}
\le 
2^{4n+11}\,\rho^{-4}\,(N^0+\bar\eps\,(t+2)^2).
\end{equation}
Now, observe that 
$$
\sum_{\tau=0}^{\sigma-1} \sum_{k=0}^{n-1} 1_{\ev X'(\tau,\sigma,k)}<3+1_{\ev V^\sigma(n)}\,n\,,
$$
since $\ev X'(\tau,\sigma,k)$ can hold for at most one $\tau$.
Therefore, by taking expectations and applying~\eref{e.fsum} we get
\begin{equation*}
\begin{aligned}
&
\sum_{\sigma=0}^{2t-1}\sum_{\tau=0}^{\sigma-1} \sum_{k=0}^{n-1} \Pb{\ev X'(\tau,\sigma,k)}
\\&\qquad\qquad
\le 6\,t+
n\sum_{\sigma=0}^{2t-1} \Ps{\ev V^\sigma(n)}
\le O(t)+
O(1)\,2^{4n}\,n\,\rho^{-4}\,(N^0+\bar\eps\,t^2).
\end{aligned}
\end{equation*}
We now assume that $2^n\le t\,\rho$.  Then the inequalities~\eref{e.evX} may be applied to
the above, giving
$$
\sum_{k=0}^{n-1}
\sum_{\tau=0}^{t-1}
\sum_{\sigma=\tau+1}^{\tau+\floor{2^{k-2}/\rho}}
 2^{-k-2}\rho^{2}\,
\Pb{y_1^\tau\ge\rho}
\le O(t)+
O(1)\,2^{4n}\,n\,\rho^{-4}\,(N^0+\bar\eps\,t^2).
$$
This implies~\eref{e.coup3}, and completes the proof.
\QED
  
\begin{corollary}\label{c.coup}
Let $\cc\in(0,1/2)$.
Let $\TT$ be a random variable with values in $\N$ which is independent from
the Markov chain $(Y^t,Z^t)$. Set $\eta:=\max\{\Ps{\TT=t}:t\in\N\}$, and
suppose that $({\bar\eps})^{1-\cc}\le\eta\le ({\bar\eps})^\cc/\max\{N^0,1\}$.
Then for all $\lambda\ge 1$ and $\rho>0$
$$
\Pb{y_1^\TT\ge\rho}
\le
\Pb{\TT>\lambda\, \eta^{-1}} + C\,(\lambda/\rho)\,\bigl|\log\bar\eps\bigr|^{-1},
$$
where $C$ is a constant depending only on $\cc$.
\end{corollary}

\proof
Let $s:=\floor{\lambda\,\eta^{-1}}$.
We have 
\begin{equation*}
\begin{aligned}
\Pb{y_1^\TT\ge\rho}
&
\le
\Pb{\TT> \lambda\,\eta^{-1}}
+
\sum_{\tau=0}^{s} \Pb{\TT=t}
\Pb{y_1^\tau\ge\rho}
\\&
\le
\Pb{\TT> \lambda \,\eta^{-1}}
+
\eta\sum_{\tau=0}^{s} 
\Pb{y_1^\tau\ge\rho}.
\end{aligned}
\end{equation*}
Thus, the proof is completed by applying~\eref{e.coup3} with
$t=s+1$ and $n:=\floor{|\log \bar\eps|/C}$,
provided that with sufficiently large $C$ we have
\begin{equation}\label{e.req}
2^{4n}\rho^{-5}\bigl(N^0/t+\bar\eps\,t\bigr)
\le
C\,\rho^{-1}\,\bigl|\log\bar\eps\bigr|^{-1}
\end{equation}
and $2^{n}\le t\,\rho$.
First, note that we may assume that $\lambda,\rho^{-1}<\bigl|\log\bar\eps\bigr|$ 
and $\bar\eps<1/10$. Then $2^{n}\le t\,\rho$ holds by the assumptions on $\eta$.
 It is also easy to verify that with an appropriate choice of $C$~\eref{e.req} follows
from our inequalities for $\eta$ and assumptions about $\rho$ and $\lambda$.
\QED
  
\begin{theorem}\label{t.mix}
Let $Y^0$ and $Z^0$ be independent random samples from $PD(1)$,
and let $(Y^t,Z^t)$ denote their evolution under $\tilde M$.
Let $t_0\in\Np$, and let $\TT\in\{0,1,\dots,t_0-1\}$ be
chosen uniformly and independently from the evolution of the
chain $(Y^t,Z^t)$.
Then for each $\rho>0$,
$$
\PB{\max\{y_1^\TT,z_1^\TT\}>\rho}\le
O(1)\,\rho^{-1}\,(\log t_0)^{-1}\,.
$$
\end{theorem}

\proof
Set $\eps:=(t_0)^{-2}$, and define $\bar\eps$ as in Lemma~\ref{l.coup1}.
Recall that a size biased sample from the $PD(1)$ sample $Y^0$ gives the
uniform distribution on $[0,1]$ (this is well-known, but also easy
to verify from the definition). 
Consequently, $\E[\bar\eps]=3\eps$.
Let $\ev A_1$ be the event that $\bar\eps\le  \eps^{3/4}$.
Then $\P[\neg \ev A_1]\le 3\,\eps^{1/4}$.
Let $\ev A_2$ be the event that $N^0\le \eps^{-1/4}$.
It is easy to see (e.g., using the
description of $PD(1)$ from the introduction)
that $\P[\neg \ev A_2]\le O(\eps)$.
(In fact, $N^0/|\log\eps|$ is very unlikely to be large.)
Define $\eta$ as in Corollary~\ref{c.coup}.
Then $\eta=\eps^{1/2}$ and on $\ev A_1\cap\ev A_2$ we have
$({\bar\eps})^{1-\cc}\le\eta\le ({\bar\eps})^\cc/\max\{N^0,1\}$
with $\cc=1/5$, for example.
On the event $\ev A_1\cap\ev A_2$,
apply the corollary with $\lambda=1$ and the corresponding statement
with the roles of $Y$ and $Z$ switched, to get
$$
\PB{\max\{y_1^\TT,z_1^\TT\}>\rho\md \ev A_1\cap\ev A_2}\le O(\rho^{-1})\,
|\log\eps|^{-1}\,.
$$
Now our estimates for $\P[\neg \ev A_1]$ and $\P[\neg\ev A_2]$
complete the proof.
\QED

\proofof{Theorem~\ref{t.DMZZ}}
The proof is similar to the proof of Theorem~\ref{t.mix}
Let $\mu$ be a measure that is invariant under $M$, and let $Y^0$
be a sample from $\mu$.  Let $Z^0$ be a sample from $PD(1)$ (which we may take
to be independent from $Y^0$, though this is not important).
Let $t_0\in\Np$, $t_0>5$, and let $\TT$ be as in Theorem~\ref{t.mix}.
As in the proof of that theorem, choose $\eps=t_0^{-2}$.

Note that for every $t\in\N$, $Y^t$ is also a sample from $\mu$,
because $\mu$ is $t$ invariant.  The same also holds for $Y^\TT$,
since $\TT$ is independent from the chain $(Y^t)$.

We now explain how to get bounds on the distributions of $\bar\eps$ and $N^0$
using continuous analogs of Lemmas~\ref{l.grow} and~\ref{l.precoup}.
Let $\beta(s,Y):=\sum\{Y_i:Y_i\le s\}$.
Since $\lim_{s\searrow0}\beta(s,Y^0)=0$ a.s.,
we may choose $k=k(\eps)>0$ sufficiently large so that
$\Pb{\beta(2^{-k},Y^0)>\eps}<\eps$ and $2^{-k}<\eps\,(1-\eps)/8$.
Set $\delta=1-2^{-k}$ and
$t_1=\ceil{2^6\,\delta^{-1}\,k\,2^k}$.
The proof of Lemma~\ref{l.grow} applied to the continuous setting gives
$$
\EB{\beta(\eps,Y^{t_1})\md \beta(2^{-k},Y^0)\le \eps}
\le O(1)\,\eps\,|\log \eps|\,.
$$
By our choice of $k$ this implies
$\Eb{\beta(\eps,Y^{t_1})}\le O(\eps)\,|\log\eps|$.
Since $Y^{t_1}$ and $Y^0$ have the same distribution, this gives
\begin{equation}\label{e.b}
\Eb{\beta(\eps,Y^{0})}\le O(\eps)\,|\log\eps|\,,
\end{equation}
and since $\Eb{\beta(\eps,Z^{0})}=\eps$, as in the proof of Theorem~\ref{t.mix}, we
conclude that $\Eb{\bar\eps}\le O(\eps)\,|\log\eps|$.

We now adapt the latter part of the proof of Lemma~\ref{l.precoup}.
Let $m_0:=\ceil{|\log_2\eps|}$.
On the one hand
$$
\sum_{m=0}^{m_0} 2^m \beta(2^{-m},Y^0)
=
\sum_{m=0}^{m_0} \sum\bigl\{ 2^m\, Y_i:Y_i\le 2^{-m}\bigr\}
\le 2\Bigl|\{i\in\Np:Y^0_i\ge \eps\}\Bigr|\,.
$$
On the other hand,~\eref{e.b} gives
$$
\EB{
\sum_{m=0}^{m_0} 2^m \beta(2^{-m},Y^0)}
\le O(1)\, \sum_{m=0}^{m_0} m = O(1)\,|\log\eps|^2\,.
$$
Consequently, we have $\E[N^0]=O(1)\,|\log\eps|^2$.
Now, the proof of Theorem~\ref{t.mix} applies, and gives for all $\rho>0$
$$
\Pb{\max\{y_1^\TT,z_1^\TT\}>\rho}\le
O(1)\,\rho^{-1}\,(\log t_0)^{-1}\,.
$$
We conclude that for every $\rho>0$ there is a coupling of $Y^0$ and $Z^0$ so that
$\Pb{\max\{y_1^0,z_1^0\}>\rho}<\rho$.  This implies that $\mu=PD(1)$.
\QED

\section{Conclusion}\label{s.conc}

\proofof{Theorem~\ref{t.main}}
The proof is similar to the proof of Theorem~\ref{t.mix}.
Let $\eps>0$. Let $\TT$ be uniformly chosen in
$(2\,\Z)\cap[0,\eps^{-1/2}]$.
Set $z=z(2\,t/n)$, as in~\eref{e.z}.
Let $Z^0$ be chosen according to $PD(1)$, and let
$Y^\tau=\SS(\pi_{t+\tau})/(n\,z)$.

We now apply a coupling of $Z^\tau$ and $Y^\tau$
similar to the coupling $\tilde M$ given in Section~\ref{s.coup}.
There are a few minor necessary modifications in the definition of the 
coupling.  
First, note that the entries of $Y^\tau$ do not sum to $1$ but to
$1/z$.
Thus, the random variables $u$ and $v$ needed in the transition
kernel for $\tilde M$ should be uniform in $[0,1/z]$.
In $\tilde Y^\tau$ and $\hat Y^\tau$ we put those segments corresponding to
cycles that do not intersect $V_G^t$ in the very end; that is, roughly in
the interval $[1,1/z]$.
When $u$ or $v$ turn out to be outside of $[0,1]$, we make no transition
to $Z$; that is, $Z^{\tau+1}=Z^\tau$, in this case.

Another modification is necessary because the transitions of $Y^\tau$ are
discrete. Thus, the actual size of the splits occuring in the transitions
of $Y$ would be determined with $\ceil{n\,z\,v}/(n\,z)$.
The definition of the matching between entries in $Y^\tau$ and entries
in $Z^\tau$ need to be modified as well. 
When a split is made in both $Y^\tau$ and $Z^\tau$, 
pieces which would have been exactly the same, may differ slightly now,
because of the discretization in the transition of $Y$. This difference
is of order $1/n$, and may be safely ignored. Although these errors
may accumulate over time, when matched pieces are merged and split,
the total discrepancy would still be small, since we take $n$ much larger
than $\eps^{-1/2}$, which bounds $\TT$.

Lemma~\ref{l.precoup} gives us good control on $N^0$
while~\eref{e.capture}
gives a bound on the probability that $\bar \eps$ is large.
Consequently, the proof of Theorem~\ref{t.mix}\ shows
that for all $\rho>0$ if $n$ is large
\begin{equation}
\label{e.nr}
\Pb{\|Y^\TT-Z^\TT\|_\infty>\rho}\le O(1)\,\rho^{-1}\,|\log\eps|^{-1}.
\end{equation}
Thus, the statement of Theorem~\ref{t.main} is obtained with $t$ replaced by $t+\TT$.
If we consider $\eps$ and $\rho$ as fixed, then $\TT$ is bounded.
Since $\TT$ is even, the following lemma completes the proof.
\QED

\begin{lemma}\label{l.fin}
Let $t\ge c\,n$, $c>1/2$.
As $n\to\infty$, the total variation distance between the law of $\SS(\pi_t)$
and the law of $\SS(\pi_{t+2})$ tends to zero.
\end{lemma}

First, we give a slightly informal proof.
Note that when the largest entry in $\SS(\pi_\tau)$ is not too small, there is probability
bounded away from $0$ and $1$ that $\SS(\pi_\tau)=\SS(\pi_{\tau+2})$, because that
entry may split and then recombine. We know that for many $\tau\in [n/2,t]$
the largest entry is not small. Consequently, there is a random \lq\lq delay\rq\rq,
which implies the statement of the lemma. For readers who are not convinced yet, we
offer a proof with more details.

\proof
Set $W^\tau=\SS(\pi_\tau)$. To prove that $W^t$ and $W^{t+2}$ have close distributions,
we couple the chain $(W^\tau)$ with a chain $(U^\tau)$ which has the same distribution
as $(W^\tau)$. In essence, the two chains will be the same; the significant difference involves
a random shift in time.
Set $\tau_0=\tau'_0=0$.
Inductively, suppose that $\tau_i$ and $\tau'_i$ have beed defined
such that $W^{\tau_i}=U^{\tau'_i}$.
Let $m=m_i$ be the largest integer such that $W^{\tau_i+2j}=W^{\tau_i}$ for all
$j=1,2,\dots,m$. The distribution of $m$ conditioned on $W^{\tau_i}$ is geometric;
that is, $\Pb{m=k\md W^{\tau_i}}=(1-p)\,p^k$, where $p=p_i=\Pb{m>0\md W^{\tau_i}}$.
Similarly, the largest integer $m'$ such that
$U^{\tau'_i+2j}=U^{\tau'_i}$ for all
$j=1,2,\dots,m'$ has the same conditioned distribution:
$\Pb{m'=k\md W^{\tau_i}}=  \Pb{m'=k\md U^{\tau'_i}}= (1-p)\,p^k$.
We now couple $m$ and $m'$.
If $\tau'_i= \tau_i+2$, take $m'=m$.
Otherwise we couple $m$ and $m'$ so that $|m-m'|\le 1$,
but $m\ne m'$ happens quite frequently.
For example, for all $k\in\N$ take 
$(m,m')=(k,k+1)$ with probability $p^{k+1}(1-p)/(1+p)$,
$(m,m')=(k+1,k)$ with the same probability,
and $(m,m')=(0,0)$ with probability $(1-p)/(1+p)$,
all conditioned on $W^{\tau_i}$.
In this case, the conditioned probability
that $m'-m=\pm 1$ is $2p/(1+p)$.
In either case, take $U^{\tau'_i+j}=W^{\tau_i+j}$ for $j=1,2,\dots,\min\{2m,2m'\}$.
If $m'>m$ let $U^{\tau'_i+2m+1}$ be independent from the chain $(W^\tau)$ given
$(W^{\tau_i},m,m')$, and similarly if $m>m'$.
Clearly, $W^{\tau_i+2m}=U^{\tau'_i+2m'}$.
Take $U^{\tau'_i+2m'+j}=W^{\tau_i+2m+j}$ for $j=1,2$,
$\tau_{i+1}:=\tau_i+2m+2$ and
$\tau'_{i+1}:=\tau'_i+2m'+2$.
Then continue inductively.
This completes the specification of the coupling.

It clearly suffices to prove that with probability tending to $1$ as
$n\to\infty$, we have $\tau'_i=\tau_i+2$ with some $\tau_i<t$.
First, observe that the $p_i$ are bounded away from $1$.
This guarantees that a.a.s.\ $\tau_i=O(i)$.
Now note that $p_i$ is bounded away from zero by some positive
function of $W^{\tau_i}_1/n$ (the largest entry normalized), because that
largest entry may split in the next step, and then the same two parts may merge
in the step after that. We know, for example from~\eref{e.capture},
that with high probability for most values 
of $i$ such that $t\ge\tau_i\ge (c\,n+n/2)/2$, the largest entry of $W^{\tau_i}$ is not too much
smaller than $n$.
Consequently, a.a.s.\ we have $p_i$ bounded away from zero for many values
of $i$ satisfying $\tau_i<t$. Similarly, $m_i\ne m'_i$ for many values of $i$.
Note that $(\tau_i-\tau'_i)/2$ is a martingale, and its increments are $\{-1,0,1\}$.
By removing the $0$ increment steps, the martingale may be coupled with a simple
random walk on $\Z$.
The martingale starts at $0$. Thus, the probability that many $\pm 1$ steps are
performed and it never gets to $1$ tends to $0$. This completes the proof.
\QED

\bigskip\noindent{\bf Acknowledgments:}
I have had the pleasure to benefit from conversations with
Rick Durrett, Michael Larsen, Russ Lyons, David Wilson and Ofer Zeitouni
in connection with this work. Nathana\"el Berestycki has kindly pointed
out an error in the proof of Lemma~\ref{l.coup3} in a previous version
of this paper.

\bibliographystyle{halpha}
\bibliography{mr,prep,notmr}

\end{document}